\documentclass{article}
\usepackage{arxiv}
\usepackage[utf8]{inputenc}
\usepackage[T1]{fontenc}
\usepackage{hyperref}     
\usepackage{booktabs}   
\usepackage{amsmath}
\usepackage{amsfonts}     
\usepackage{nicefrac}     
\usepackage{microtype}
\usepackage{lipsum}
\usepackage{bm}
\usepackage{amsmath}
\usepackage{multirow}
\usepackage{graphicx}
\graphicspath{ {./images/} }
\usepackage{subcaption}
\usepackage[authoryear,round]{natbib}
\usepackage{xcolor}
\usepackage{soul}
\RequirePackage[thmmarks]{ntheorem}

\title{Transformed-Linear Innovations Algorithm for Modeling and Forecasting of Time Series Extremes}

\author{
  Nehali Mhatre\\
  Department of Statistics\\
  Colorado State University\\
  \texttt{nehalimhtr@gmail.com} \\
   \And
 Daniel Cooley \\
  Department of Statistics\\
  Colorado State University\\
  \texttt{cooleyd@stat.colostate.edu} \\
}

\theoremheaderfont{\hskip\parindent\normalfont\scshape}
\theorembodyfont{\itshape}
\theoremseparator{.}
\newtheorem{proposition}{Proposition}
\newtheorem{theorem}{Theorem}
\newtheorem{remark}{Remark}
\newtheorem{lemma}{Lemma}
\newtheorem{corollary}{Corollary}

\begin{document}
\maketitle




\begin{abstract}
The innovations algorithm is a classical recursive forecasting algorithm used in time series analysis.
We develop the innovations algorithm for a class of nonnegative regularly varying time series models constructed via transformed-linear arithmetic.
In addition to providing the best linear predictor, the algorithm also enables us to estimate parameters of transformed-linear regularly-varying moving average (MA) models, thus providing a tool for modeling.

We first construct an inner product space of transformed-linear combinations of nonnegative regularly-varying random variables and prove its link to a Hilbert space which allows us to employ the projection theorem, from which we develop the transformed-linear innovations algorithm.
Turning our attention to the class of transformed linear MA($\infty$) models, we give results on parameter estimation and also show that this class of models is dense in the class of possible tail pairwise dependence functions (TPDFs).
We also develop an extremes analogue of the classical Wold decomposition.
Simulation study shows that our class of models captures tail dependence for the GARCH(1,1) model and a Markov time series model, both of which are outside our class of models.

We also develop prediction intervals based on the geometry of regular variation.
Simulation study shows that we obtain good coverage rates for prediction errors.
We perform modeling and prediction for hourly windspeed data by applying the innovations algorithm to the estimated TPDF.

\end{abstract}



\section{Introduction}
\label{sec:introduction}
A primary aim of time series analysis is forecasting.
\citet{mhatrecooley2020} developed transformed-linear time series models, a class of time series which are nonnegative and regularly-varying and which are similar to familiar ARMA models in the non-extreme setting.
\citet{mhatrecooley2020} showed that these relatively simple models can capture dependence in a time series' upper tail. 
We now address the problem of forecasting specifically when values are large.

Our approach for forecasting is to develop the innovations algorithm for transformed-linear time series.
The innovations algorithm, a well known forecasting method in classical time series, relies on the autocovariance function.
In traditional time series analysis and elsewhere, the best linear predictor $\hat{X}_{n+1}$ minimizes mean squared prediction error (MSPE), $E[(X_{n+1} - \hat{X}_{n+1})^2]$ and Gaussian assumptions are usually used to create  prediction intervals.
However, the autocovariance function is not well-suited for describing tail dependence and expected-squared error is not a natural or intuitive measure of loss for extremes.
The innovations algorithm we develop uses the tail pairwise dependence function (TPDF) to characterize extremal dependence, and minimizes a quantity describing tail behavior.
Despite these differences, we show that the form of the transformed-linear predictor is of the same form as in the non-extreme setting.

To develop the innovations algorithm, we construct a vector space $\mathbb{V}$ of a series of absolutely summable transformed-linear combinations of nonnegative regularly-varying random variables.
We show that $\mathbb{V}$ is an inner product space and is isomorphic to $\ell^1$, the space of absolutely summable sequences.
Although $\mathbb{V}$ itself is not a Hilbert space, we show that the set of predictors based on previous $n$ observations is isomorphic to a closed linear subspace of $\ell^2$, the space of square summable sequences, and we can employ the projection theorem.
Using the properties of the projection theorem we develop a transformed-linear analogue of the classical innovations algorithm that allows us to do modeling and prediction iteratively.

The innovations algorithm gives us more than just a method for prediction, as we also use it as a tool to demonstrate properties of this modeling framework.
Using the innovations algorithm we show that if the true model is in our transformed-linear space then applying the innovations algorithm iteratively yields parameter estimates that converge to the true parameters.
Furthermore we show that even if the underlying model is not a transformed-linear model, applying the innovations algorithm will yield a transformed-linear model whose TPDF matches closely the estimated TPDF of the underlying model.
We go on to show that the class of transformed-linear regularly-varying MA($\infty$) time series is dense in the class of possible TPDFs.
We also develop a transformed-linear analogue of the Wold decomposition.
To demonstrate the richness of the class of transformed-linear regularly-varying MA($\infty$) models we run the innovations algorithm on data simulated from two different models from outside our class of transform linear models: the GARCH(1,1) process and a first-order Markov chain.
Neither of these models is in the family of transformed-linear time series.
We show for both these models that by running the innovations algorithm on the estimated TPDF we can get estimates for coefficients of a transformed-linear regularly-varying MA time series whose TPDF closely matches the estimated TPDF of the simulated data and well represents summary measures of tail dependence.

We also develop a method based on the polar geometry of regular variation for producing prediction intervals for the case when predictands are large.
Because the regular variation geometry differs from the elliptical geometry typically assumed in standard linear prediction settings, uncertainty quantification is significantly different from the non-extreme setting.
We perform modeling and prediction for the windspeed anomalies data discussed in \citet{mhatrecooley2020} by applying the innovations algorithm to the estimated TPDF.


\section{Background: Transformed-Linear Time Series Models for Extremes}
\label{sec:background}

\cite{mhatrecooley2020} use transformed-linear arithmetic to produce nonnegative regularly varying time series models for capturing extremal dependence.  
Here we review the basics of these models.
We consider time series $\{X_t\}, t \in \mathbb(Z)$ whose finite-dimensional distributions are multivariate regularly varying.
Let ${\bm X}_{t,p} = (X_t, X_{t+1}, \ldots, X_{t+p-1})^T$ for any $t$ and $p > 0$.
Then, there exists a function $b(s) \rightarrow \infty$ as $s \rightarrow \infty$ and a non-trivial limit measure $\nu_{\bm X_{t,p}}$ such that
\begin{align}
  s \text{ pr}\left\{ \frac{{\bm X_{t,p}}}{b(s)} \in \cdot\right\} \xrightarrow{v} \nu_{\bm X_{t,p}}(\cdot) \mbox{ as } s \rightarrow \infty,
  \label{eq:mvRegVar}
\end{align}
where $\xrightarrow{v}$ denotes vague convergence in $M_+(\mathbb{R}^p \setminus \{ 0 \})$, the space of nonnegative Radon measures on $\mathbb{R}^p \setminus \{ 0 \}$ \citep[Section 6]{resnick2007}.
The normalizing function $b(s) = s^{1/\alpha} U(s)$ for $\alpha > 0$ and some slowly varying function $U$.
The scaling property $\nu_{\bm X_{t,p}}(aC) = a^{-\alpha}\nu_{\bm X_{t,p}}(C)$ for any $a > 0$ and any set $C\subset \mathbb{R}^p \setminus \{ 0 \}$ implies that the measure is more easily understood via polar coordinates.  
In particular, for a radially-defined set $C(r,B) = \{{\bm X_{t,p}} \in \mathbb{R}^p : \|{\bm X_{t,p}}\| >r, \|{\bm X_{t,p}}\|^{-1}{\bm X_{t,p}} \in B\}$ where Borel set $B \subset \mathbb{S}^{p-1}= \{\bm x \in \mathbb{R}^p : \|\bm x\|=1\}$, $\nu_{\bm X_{t,p}}(C(r,B)) = r^{-\alpha}H_{\bm X_{t,p}}(B)$ where $H_{\bm X_{t,p}}$ is the angular measure taking values on $\mathbb{S}^{p-1}$.

\citet{kuliksoulier2020} provide a comprehensive treatment of regularly-varying time series, beginning with the finite-dimensional distribution definition and developing the limit measure of the process.
Fully characterizing the limit measure is challenging, particularly because only a subset of extreme data are used for estimation.
\cite{mhatrecooley2020} instead develop the notion of tail stationarity, an extremal analogue to second-order stationarity.
Characterization of the tail of a time series is simplified to characterizing the pairwise dependence via the TPDF, which takes the place of the autocovariance function in standard time series analysis.  
The TPDF is the functional extension of the tail pairwise dependence matrix (TPDM), introduced in \cite{cooleythibaud2019}.
\cite{kiriliouk2022} define the TPDM for regularly-varying random vectors with tail index $\alpha$, but we will follow both \cite{cooleythibaud2019} and \cite{mhatrecooley2020} and further assume $\{X_t\}$ has tail index $\alpha = 2$, so that in Section \ref{sec:subsetVplus} we can connect the TPDF to the inner product we introduce in Section \ref{sec:vectorspace}.
In practice, this assumption generally requires a marginal transformation.
The TPDF is given by
\begin{eqnarray} 
\label{eq:tpdf}
\sigma(X_t, X_{t+h}) = \int_{\Theta_{1}} w_t w_{t+h}\text{d}H_{X_t, X_{t+h}}(\bm w), 
\end{eqnarray}
where $\Theta_{1} = \{\bm x \in \mathbb{R}^2 \mid \| \bm x \|_2 = 1\}$.
If the time series is stationary the TPDF can be viewed as a function of the lag $h$.
Like the autocovariance function in standard time series analysis, the TPDF evaluated at lag 0 provides a measure of the marginal `scale' of the time series, as
$
\lim \limits_{s\to \infty}  s \text{ pr} \left\{ \frac{|X_1|}{b(s)} > c \right\} 
  = c^{-2} \sigma(0).
$
As $\sigma(h)/\sigma(0) \in [0,1]$, this ratio provides an interpretable number for dependence strength at lag $h$.


Linear time series models of the form $X_t = \sum_{j = -\infty}^\infty \psi_j Z_{t-j},$ where $Z_{t-j}$ is a white noise sequence make up a large portion of classical time series analysis, and include the familiar ARMA models.
For time series with finite second moments, $n$-step predictors (predictors of $X_{t+1}$ based on $X_t, \ldots X_{t-n+1})$ minimizing mean square prediction error (MSPE) can be constructed via the projection theorem.
It is straightforward to construct linear regularly varying time series models by assuming the noise sequence is regularly varying.
As in \cite{mhatrecooley2020}, we are motivated to construct {\em nonnegative} time series models, which enables one to focus attention on the time series upper tail.
If one wants to construct a nonnegative linear time series using traditional arithmetic, one needs to restrict $\psi_j > 0$ for all $j$ and employ a nonnegative noise sequence.
\cite{mhatrecooley2020} show that transformed-linear time series can more flexibly fit data than traditional linear time series models restricted to be nonnegative.
Important for this work, allowing the $\psi_j$'s to take negative values will allow us to show the elements of transformed-linear $\{X_t\}$ can be thought of as members of a vector space.

\cite{mhatrecooley2020} employ the transformed-linear arithmetic of \cite{cooleythibaud2019} to construct nonnegative time series models.
Given bijective transform $\tau: \mathbb{R} \mapsto \mathbb{R}_+$,  
define $X_1 \oplus X_2 = \tau(\tau^{-1}(X_1) + \tau^{-1}(X_2))$ and $a \circ X_1 = \tau(a \tau^{-1}(X_1)$ for $a \in \mathbb{R}$.
\cite{cooleythibaud2019} show that if $X_1, X_2$ are independent regularly varying with index $\alpha$ with respective measures $\nu_{X_1}, \nu_{X_2}$, if $\lim_{y \rightarrow \infty} \tau(y)/y = 1$, and if $X_1, X_2$ meet a lower tail condition associated with the particular $\tau$, then $\nu_{X_1 \oplus X_2}(\cdot) = \nu_{X_1}(\cdot) + \nu_{X_2}(\cdot)$ and $\nu_{a \circ X_1}(\cdot) = (a^{(0)})^\alpha \nu_{X_1}(\cdot)$, where $a^{(0)} = \max(a, 0)$.
If $\tau(y) = \log{\{1 + \exp{(y)}\}}$, the lower tail condition
$
s\text{ pr}\left[X_i \le \exp{\{-kb(s)\}}\right] \to 0, \; k>0, \; i = 1,2 \; s \rightarrow \infty
$
is sufficient.

\cite{mhatrecooley2020} largely focus on causal transformed linear time series models  
\begin{eqnarray} 
\label{eq:transLinearTS}
X_t = \bigoplus_{j=0}^{\infty}\psi_j \circ Z_{t-j},
\end{eqnarray}
where $\sum^{\infty}_{j=-\infty}|\psi_j|<\infty$, and $\{Z_t\}$ is a noise sequence of independent and tail stationary $RV_+(2)$ random variables. 
The focus of \cite{mhatrecooley2020} is to develop transformed linear ARMA models and to develop their properties, particularly the TPDF,  as this is the critical parameter from the tail stationarity standpoint.
The model in (\ref{eq:transLinearTS}) is known as the transformed-linear MA($\infty$) and has TPDF $\sigma(h) = \sum^{\infty}_{j=0} \psi^{(0)}_j \psi^{(0)}_{j+h}$.


\section{Space $\mathbb{V}$ and Innovations Algorithm}
\label{sec:VandInnovations}

\subsection{Inner Product Space $\mathbb{V}$}
\label{sec:vectorspace}

We begin by considering the space $\mathbb{V} = \{X_t : X_t = \bigoplus^{\infty}_{j=0}\psi_{t,j} \circ Z_j, \sum^{\infty}_{j=0} |\psi_{t,j}|<\infty $
where $Z_j$'s are independent and tail stationary $RV_+(2)$ random variables with $\lim_{x \to \infty} \text{Pr}(Z_j > x)/\{x^{-2}L(x)\} = 1$ for some slowly-varying function $L(x)$, $\psi_{t,j}\in \mathbb{R}$, and $t\in \mathbb{Z}$.
\cite{mhatrecooley2020} show that $X_t$ converges with probability one.
Letting $X_t = \bigoplus^{\infty}_{j=0}\psi_{t,j} \circ Z_{j}, X_s = \bigoplus^{\infty}_{j=0}\psi_{s,j} \circ Z_{j} \in \mathbb{V}$, then $X_t \oplus X_s = \bigoplus^{\infty}_{j=0}(\psi_{t,j} + \psi_{s,j}) \circ Z_{j}$ and $a \circ X_t = \bigoplus^{\infty}_{j=0} (a^{(0)} \psi_{t,j}) \circ Z_{j}$.
The space $\mathbb{V}$ equipped with transformed linear operations is a vector space.
The details of this and other properties of $\mathbb{V}$ described in this section are provided in the supplementary materials.

Define the inner product between $X_t$ and $X_{s}$ as,
\begin{eqnarray} \label{eq:innerprod}
\langle X_t, X_s \rangle := \sum^{\infty}_{j=0} \psi_{t,j}\psi_{s,j}.
\end{eqnarray}
We define the norm $\|X_t\| = \sqrt{\langle X_t, X_t \rangle}=\sqrt {\sum^{\infty}_{j=0} \psi_{t,j}^2}$ and say $X_t$ and $X_s$ are orthogonal if $\langle X_t, X_s \rangle = 0$.
Note the space $\mathbb{V}$ is more general than the stationary time series construction given in (\ref{eq:transLinearTS}).
We will return to the stationary time series setting in Section \ref{sec:statyTS}, and we will relate the inner product to the TPDF in Section \ref{sec:subsetVplus}.

It will be useful to simultaneously consider the infinite-dimensional space of absolutely summable sequences,
$$\ell^1 = \left\{\{ a_j\}^{\infty}_{j=0}, a_j \in \mathbb{R} : \sum^{\infty}_{j=0}|a_j|< \infty \right\},$$
equipped with the standard vector addition and scalar multiplication. 
For any $X_t \in \mathbb{V}$ we can define a mapping $T:\mathbb{V} \to \ell^1$ such that $T(X_t) = \{\psi_{t,j}\}_{j=0}^{\infty} \in \ell^1$.
The mapping $T$ is a linear map and an isomorphism.
We know that vector space $\ell^1 \subset \ell^2$, where $\ell^2=\{\{d_j\}: \sum^{\infty}_{j=0} d_j^2<\infty\}$, the space of square-summable sequences.
The inner product defined in (\ref{eq:innerprod}) is isomorphic to the usual inner product on $\ell^2$.
\subsection{Best Transformed-Linear Prediction in $\mathbb{V}$}
\label{sec:prediction}

We want to use the projection theorem to perform prediction.
Unfortunately, $\mathbb{V}$ is not itself a Hilbert space as $\mathbb{V}$ is isomorphic to $\ell^1$, and $\ell^1$ is not complete in the metric induced by the $\ell^2$ inner product in (\ref{eq:innerprod}).
We consider the sequence $\{X_t\}, t \in \mathbb{Z}$ and consider transformed-linear prediction in terms of the previous $n$ steps;
that is we consider predictors
\begin{eqnarray}\label{eq:predictor}
\hat{X}_{n+1}(\bm b_n) = \bigoplus^{n}_{j=1}b_{nj} \circ X_{n+1-j},
\end{eqnarray}
where $\boldsymbol{b}_n = (b_{n1}, \cdots, b_{nn})^T \in \mathbb{R}^n.$
Let $\mathbb{V}^n$ be the set of all such predictors.
Let us consider the analogous problem in $\ell^1$.
Consider sequences $\{a_j\} = T(X_{n+1-j}), j = 1, \ldots, n$ in $\ell^1$. 
Let $\mathbb{C}^n$ be the set of sequences $c(\bm b_n) = b_{n1} \{a_1\} + \cdots + b_{nn} \{a_n\}.$
$\mathbb{C}^n$ is the space spanned by $\{a_1\}, \cdots, \{a_n\}$, and dim($\mathbb{C}^n)\le n$.
Since any $n$-dimensional subspace of a complex topological vector space is closed \citep[Theorem 1.21]{rudin1991}, 
$\mathbb{C}^n$ is a closed subspace of $\ell^1 \subset \ell^2$.
By the projection theorem, there is a unique $\hat{c} \in \mathbb{C}^n$ such that $\|x-\hat{c}\| = \inf_{c\in \mathbb{C}^n}\|x-c\|$, for every $x$ in $\ell^2$.
Thus, the set of predictors $\mathbb{V}^n$ based on previous $n$ observations is isomorphic to a closed linear subspace of $\ell^2$ and we can employ the projection theorem since $\ell^2$ is known to be a Hilbert space.

Armed now with the projection theorem, the best linear one-step predictor $\hat{X}_{n+1}$ is given by
\begin{equation}
  \hat{X}_{n+1} = \left\{
                    \begin{array}{l l}
                    0, & \text{ if } n=0, \\
                    P_{\mathbb{V}^n}X_{n+1} & \text{ if } n\ge 1,
                    \end{array}
                  \right.
\end{equation}
where 
$P_{\mathbb{V}^n}$ denotes the projection mapping onto $\mathbb{V}^n$.
Thus, $\hat{X}_{n+1}$ is a transformed-linear combination of $\{X_1, ..., X_n\}$ as given in (\ref{eq:predictor}).
Define $\hat {\bm b}_n = (\hat b_{n1}, \ldots, \hat b_{nn})^T$ to be the solutions to the prediction equations given by the projection theorem
\begin{align} \label{eq:oneStep}
    \left \langle X_{n+1} \ominus \bigoplus^n_{j=1}\hat b_{nj} \circ X_{n+1-j}, \; \;  X_{n+1-k} \right \rangle & =0, \; \; k=1, \cdots, n.
\end{align}
Equivalently,
\begin{align}
    \label{eq:predEqn}
    \left \langle \bigoplus^n_{j=1}\hat b_{nj} \circ X_{n+1-j}, \; \;  X_{n+1-k} \right \rangle &= \langle X_{n+1}, \; X_{n+1-k} \rangle , \; \; k=1, \cdots, n.
\end{align}
By linearity of the inner product, the prediction equations can be rewritten in matrix form as
\begin{align}
    \label{eq:oneStepMatrix}
    \Gamma_n \hat {\bm b}_n & = \boldsymbol{\gamma}_n
\end{align}
where $\Gamma_n = \left[ \langle X_{n+1-j}, X_{n+1-k} \rangle \right]^n_{j,k=1}$, 
and $\boldsymbol{\gamma}_n = \left[ \langle X_{n+1}, X_{n+1-k} \rangle \right]^n_{k=1}$. If $\Gamma_n$ is non-singular, then the solution is given as
\begin{align}
    \label{eq:bn}
    \hat{\boldsymbol{b}}_n & = \Gamma_n^{-1}\boldsymbol{\gamma}_n.
\end{align}
It can be shown that the above is equivalent to minimizing the squared norm $\|X_{n+1} \ominus \hat{X}_{n+1} \|^2$ by setting the appropriate derivative to zero.
We see that (\ref{eq:bn}) is of the familiar form for linear prediction in the non-extreme setting where the inner product terms are autocovariances.

\subsection{Transformed-Linear Innovations} \label{sec:innovations}

Following \citet{brockwelldavis1991}, we develop a transformed-linear analogue of the recursive innovations algorithm to obtain the one-step predictors $\hat{X}_{n+1}$, $n\ge1$, defined in (\ref{eq:oneStep}), without having to perform matrix inversion of $\Gamma_n$.

Consider the transformed-linear innovation, $(X_{n+1} \ominus \hat{X}_{n+1})$, $n \ge 1$.
Since $\mathbb{V}^n = \Bar{\text{sp}}\{X_1, \cdots, X_n\}$, letting $\hat{X_1}:=\tau^{-1}(0)$, $\mathbb{V}^n=\Bar{\text{sp}}\{X_1 \ominus \hat{X_1}, \cdots, X_n \ominus \hat{X}_n\}$.
We can rewrite the predictor in (\ref{eq:predictor}) in terms of the innovations as,
\begin{eqnarray}
  \hat{X}_{n+1} = \bigoplus_{j=1}^n \theta_{nj} \circ \left(X_{n+1-j} \ominus \hat{X}_{n+1-j}\right). \nonumber
\end{eqnarray}
By properties of projection mappings, $\hat{X}_{n+1} \in \mathbb{V}^n$ and by (\ref{eq:predEqn}),
\begin{eqnarray}
  \langle X_{n+1} \ominus \hat{X}_{n+1},  \hat{X}_{n+1} \rangle = 0. \nonumber
\end{eqnarray}
That is, the transformed-linear innovation $(X_{n+1} \ominus \hat{X}_{n+1})$ is orthogonal to a transformed-linear combination of $X_1, ..., X_n$.
Thus, the innovation is orthogonal to each of $X_1, ..., X_n$.

Consider the set of transformed-linear innovations, $\{X_{n+1-j} \ominus \hat{X}_{n+1-j}\}_{j=1, ..., n}$. The innovation $(X_i \ominus \hat{X}_i) \in \mathbb{V}^{j-1}$ for $i < j$, as $(X_i \ominus \hat{X}_i)$ is a transformed-linear combination of $X_1, ..., X_i$.
Also, by (\ref{eq:predEqn}), $(X_j \ominus \hat{X}_j) \perp \mathbb{V}^{j-1}$.
Thus, elements of the set $\{X_1 \ominus \hat{X}_1, X_2 \ominus \hat{X}_2, ..., X_n \ominus \hat{X}_n\}$ are mutually orthogonal.
In fact, $\{X_{n+1-j} \ominus \hat{X}_{n+1-j}\}_{j=1, ..., n}$ is an orthogonal basis of $\mathbb{V}^n$.
Let the squared distance of prediction be denoted by $\nu_n$, that is, $\nu_n = \|X_{n+1} \ominus \hat{X}_{n+1}\|^2$.
The innovations algorithm  for a transformed-linear time series in $\mathbb{V}$ is given as follows:

\begin{proposition}[The Transformed-Linear Innovations Algorithm]
If $\{X_t\}$ is a transformed-linear time series in $\mathbb{V}$, where the matrix $\Gamma_n = \left[ \langle X_i,X_j \rangle \right]^n_{i,j=1}$ is non-singular for each $n \ge 1$, then the one-step predictors $\hat{X}_{n+1}$, $n \ge 0$, and their squared distances of prediction $\nu_n$, $n \ge 1$, are given by
\begin{eqnarray} \label{x_hat}
  \hat{X}_{n+1} =
    \begin{cases}
      0 & \text{if } n=0\\
      \bigoplus_{j=1}^n \theta_{nj} \circ (X_{n+1-j} \ominus \hat{X}_{n+1-j}) & \text{if } n \ge 1,
    \end{cases}       
\end{eqnarray}
and
\begin{eqnarray} \label{ia}
    \begin{cases}
      \nu_0 & = \langle X_1,X_1 \rangle\\
      \theta_{n,n-k} & = \nu_k^{-1}\Big(\langle X_{n+1}, X_{k+1} \rangle - \sum_{j=0}^{k-1}\theta_{k, k-j} \theta_{n, n-j} \nu_j \Big), \quad k=0,1,...,n-1,\\
      \nu_n & = \langle X_{n+1}, X_{n+1} \rangle - \sum_{j=0}^{n-1} \theta_{n,n-j}^2 \nu_j,
    \end{cases}.       
\end{eqnarray}
\end{proposition}

\noindent{Proof}: 
Taking the inner product on both sides of (\ref{x_hat}) with $(X_{k+1} \ominus \hat{X}_{k+1})$, $0 \le k < n$, we get
\begin{align*}
 \left\langle \hat{X}_{n+1}, (X_{k+1} \ominus \hat{X}_{k+1}) \right\rangle & = \left\langle\left\{ \bigoplus_{j=1}^n \theta_{nj} \circ (X_{n+1-j} \ominus \hat{X}_{n+1-j}) \right\},  (X_{k+1} \ominus \hat{X}_{k+1}) \right\rangle \\
 & = \sum_{j=1}^n \theta_{nj} \left\langle(X_{n+1-j} \ominus \hat{X}_{n+1-j}),  (X_{k+1} \ominus \hat{X}_{k+1}) \right\rangle \\
 & = \theta_{n, n-k} \nu_k,
\end{align*}
since $(X_{n+1-j} \ominus \hat{X}_{n+1-j}) \perp (X_{k+1} \ominus \hat{X}_{k+1})$ for all $j \ne n-k$.

Also, since $(X_{n+1} \ominus \hat{X}_{n+1}) \perp (X_{k+1} \ominus \hat{X}_{k+1})$ for $k=0, \cdots, n-1$, we get, $$\langle \hat{X}_{n+1}, (X_{k+1} \ominus \hat{X}_{k+1}) \rangle = \langle X_{n+1}, (X_{k+1} \ominus \hat{X}_{k+1}) \rangle.$$
Hence, the coefficients $\theta_{n, n-k}$, $k=0, ..., n-1$ are given by
\begin{eqnarray} \label{thetas}
 \theta_{n,n-k} = \nu_k^{-1} \langle X_{n+1}, (X_{k+1} \ominus \hat{X}_{k+1}) \rangle.        
\end{eqnarray}
Using the representation in (\ref{x_hat}) with $n$ replaced by $k$, we get
\begin{eqnarray} \label{thetas_2}
 \theta_{n,n-k} = \nu_k^{-1} \Big( \langle X_{n+1}, X_{k+1} \rangle - \sum_{j=0}^{k-1} \theta_{k, k-j} \langle X_{n+1}, (X_{j+1} \ominus \hat{X}_{j+1}) \rangle \Big).      
\end{eqnarray}
Since by (\ref{thetas}), $\langle X_{n+1}, (X_{j+1} \ominus \hat{X}_{j+1}) \rangle = \nu_j \theta_{n, n-j}$, $0 \le j < n$, we can rewrite (\ref{thetas_2}) as
\begin{eqnarray}
 \theta_{n,n-k} = \nu_k^{-1} \Big( \langle X_{n+1}, X_{k+1}\rangle - \sum_{j=0}^{k-1} \theta_{k, k-j} \theta_{n, n-j} \nu_j \Big). \nonumber
 \end{eqnarray}
By properties of projection mapping, we have
\begin{eqnarray}
 \nu_n = \|X_{n+1} \ominus \hat{X}_{n+1}\|^2 = \|X_{n+1}\|^2 - \|\hat{X}_{n+1}\|^2 = \langle X_{n+1},X_{n+1}\rangle - \sum_{k=0}^{n-1} \theta_{n,n-k}^2 \nu_k. \nonumber
\end{eqnarray}
\hspace*{0pt}\hfill $\square$

We will return to forecasting in Section \ref{sec:pred_error} where we will propose a method for quantifying prediction uncertainty.
For now, we will turn our attention to using the innovations algorithm as a tool to better understand the richness of our class of models and also as a tool for model fitting.

\section{Implications for Modeling of Stationary Time Series}
\label{sec:statyTS}


Rather than the general setting described by the class $\mathbb{V}$, we now focus on stationary time series.
If $\{X_t\}$ is an MA($\infty$) time series (\ref{eq:transLinearTS}), $X_t \in \mathbb{V}$ for all $t$.
As $\{X_t\}$ is stationary, it is natural to think of the inner product as a function of lag:
$$
\gamma(h) = \langle X_t, X_{t + h} \rangle = \sum^{\infty}_{j=0}\psi_j \psi_{j+h}.
$$
Being an inner product, $\gamma(.)$ is nonnegative definite and by the Cauchy-Schwarz inequality, $|\gamma(h)|\leq \gamma(0)$. 


The following corollary shows that given an invertible transformed-linear regularly-varying MA time series, iterating the transformed-linear innovations algorithm yields coefficient estimates that converge to the true MA coefficients.

\begin{corollary}\label{prop:convergence}
If $\{X_t\}$ is an invertible MA process in $\mathbb{V}$, that is, $$Z_t = X_t \oplus \bigoplus^{\infty}_{j=1} \pi_j \circ X_{t-j},$$
with $\lim_{x \to \infty} \text{Pr}(Z_j > x)/\{x^{-2}L(x)\} = 1$, then as $n\to \infty$, \\
(i) $\nu_n \to 1$, \\
(ii) $\|(X_n \ominus \hat{X}_n) \ominus Z_n\|^2 \to 0$, and \\
(iii) $\theta_{nj} \to \psi_j, j=1,2,\cdots$.
\end{corollary}

\noindent{Proof}:\\
Let $\mathbb{M}_{n} = \Bar{\text{sp}}\{X_s, -\infty < s\le n\}$ and $\mathbb{V}^{n} = \Bar{\text{sp}}\{X_1, \cdots, X_{n}\}$.
Because $\{X_t\}$ is invertible,
$$
Z_{n+1} \ominus X_{n+1} = \bigoplus^{\infty}_{j=1} \pi_j \circ X_{n+1-j} =  P_{\mathbb{M}_{n}}(Z_{n+1} \ominus X_{n+1}) = \ominus P_{\mathbb{M}_n}X_{n+1}, \nonumber  
$$
since $Z_{n+1} \perp \mathbb{M}_n.$
Also, we can think of $Z_k$ as $Z_k=\bigoplus^{\infty}_{j=0} \psi_j \circ Z_j$, where $\psi_j = 1$ for $j=k$ and $\psi_j=0$ for all $j\ne k$. Thus, $Z_k\in \mathbb{V}$ and subsequently, $\|Z_k\|^2 = 1$ for all $k$.
Then,
\begin{eqnarray}
1 = \|Z_{n+1}\|^2 
& =& \|X_{n+1} \oplus \bigoplus^{\infty}_{j=1} \pi_j \circ X_{n+1-j} \|^2 = \|X_{n+1} \ominus P_{\mathbb{M}_n}X_{n+1}\|^2 \nonumber \\
& \le& \|X_{n+1} \ominus P_{\mathbb{V}^n}X_{n+1}\|^2 = \nu_{n} \nonumber \\
& \le& \|X_{n+1} \oplus \bigoplus^{n}_{j=1} \pi_j \circ X_{n+1-j}\|^2 = \|Z_{n+1} \ominus \bigoplus^{\infty}_{j=n+1} \pi_j \circ X_{n+1-j}\|^2 \nonumber \\
& =& \|Z_{n+1}\|^2 + \|\bigoplus^{\infty}_{j=n+1} \pi_j \circ X_{n+1-j}\|^2 = 1 + \sum^{\infty}_{i,j=n}\pi_i\pi_j \langle X_{n+1-i}, X_{n+1-j} \rangle \nonumber \\
& \le& 1+ \left( \sum^{\infty}_{j=n+1} \pi_j \right)^2 \gamma(0). \nonumber
\end{eqnarray}
Thus (i) is established since, $$1 \le \nu_n \le 1+\left( \sum^{\infty}_{j=n+1} \pi_j \right)^2 \gamma(0) \implies \nu_n \to 1 \text{ as } n\to \infty.$$

Consider,
\begin{eqnarray} \label{proof(b)}
    \|X_n \ominus \hat{X}_n \ominus Z_n\|^2 
    & =& \|X_n \ominus \hat{X}_n\|^2 - 2\langle Z_n, X_n \ominus \hat{X}_n \rangle + \|Z_n\|^2 \nonumber \\
    & =& v_{n-1} - 2\left [\langle Z_n, X_n \rangle - \langle Z_n, \hat{X}_n \rangle \right ] +1 \nonumber \\
    & =& v_{n-1} - 2\left [\langle Z_n, \bigoplus^{\infty}_{j=0} \psi_j \circ Z_{n-j} \rangle - \langle Z_n, \bigoplus^{n-1}_{j=1}b_{nj} \circ X_{n-j} \rangle \right ] +1 \nonumber \\
    & =& v_{n-1} + 2[\|Z_n\|^2 - 0] +1 \nonumber \\
    & =& v_{n-1} - 1,
\end{eqnarray}
where (\ref{proof(b)}) converges to $0 \text{ as } n\to \infty$ by (i), thus proving (ii).

Since $X_{n+1} = \bigoplus^{\infty}_{j=0}\psi_j \circ Z_{n+1-j}$, we have that $$\psi_j = \langle X_{n+1}, Z_{n+1-j} \rangle.$$ Also, by (\ref{thetas}), $$\theta_{nj} = \nu_{n-j}^{-1} \langle X_{n+1}, (X_{n+1-j} \ominus \hat{X}_{n+1-j}) \rangle.$$
Then,
\begin{eqnarray}
    |\theta_{nj}-\psi_j| 
    & =& \left |\theta_{nj} - \langle X_{n+1}, (X_{n+1-j} \ominus \hat{X}_{n+1-j}) \rangle +\langle X_{n+1}, (X_{n+1-j} \ominus \hat{X}_{n+1-j}) \rangle - \psi_j\right | \nonumber \\
    & \le& \left |\theta_{nj} - \langle X_{n+1}, (X_{n+1-j} \ominus \hat{X}_{n+1-j}) \rangle \right | + \left |\langle X_{n+1}, (X_{n+1-j} \ominus \hat{X}_{n+1-j}) \rangle - \langle X_{n+1}, Z_{n+1-j} \rangle\right | \label{proof(c)1} \\
    & =& \left |\theta_{nj} - \theta_{nj}v_{n-j} \right | + \left |\langle X_{n+1}, (X_{n+1-j} \ominus \hat{X}_{n+1-j} \ominus Z_{n+1-j}) \rangle\right | \nonumber \\
    & \le& \left |\theta_{nj} - \theta_{nj}v_{n-j} \right | +  \sqrt{\gamma(0)}\left \|(X_{n+1-j} \ominus \hat{X}_{n+1-j} \ominus Z_{n+1-j}) \right \|, \label{proof(c)2}
\end{eqnarray}
where the inequalities in (\ref{proof(c)1}) and (\ref{proof(c)2}) hold by the triangle inequality and the Cauchy-Schwarz inequality, respectively.
Since $\theta_{nj}$ and $\gamma(0)$ are bounded, as $n\to \infty$, the first term on the right-hand side of (\ref{proof(c)2}) converges to $0$ by (i) and the second term on the right-hand side of (\ref{proof(c)2}) converges to $0$ by (ii). Thus, $\theta_{nj} \to \psi_j$ as $n \to \infty$, proving (iii). \hfill{$\square$}

Analogous to Proposition 3.2.1 in \citet{brockwelldavis1991}, Proposition \ref{prop:qcorrelated} shows that a $q$-tail-dependent tail stationary regularly-varying time series can be represented as a transformed-linear regularly-varying MA($q$) process.  
This proposition is a beginning step in showing the richness of the class of transformed-linear time series, which we build on in Section \ref{denseness}.

\begin{proposition}
  \label{prop:qcorrelated}
If $\{X_t\}$ is a regularly-varying tail stationary process in $\mathbb{V}$ with inner product function $\gamma(.)$ such that $\gamma(h)=0$ for $|h|>q$ and $\gamma(q)\ne 0$, then $\{X_t\}$ is a transformed-linear regularly-varying MA($q$) process, i.e. there exists a regularly-varying noise sequence $\{Z_t\}$ of independent and tail stationary $Z_t$'s such that
$$
    X_t = Z_t \oplus \theta_1 \circ Z_{t-1} \oplus \cdots \oplus \theta_q \circ Z_{t-q}. \nonumber
$$
\end{proposition}

Before we prove Proposition \ref{prop:qcorrelated}, we prove the following Lemma which requires the notion of convergence in tail ratio.
Let $X_t$ be the MA($\infty$) in (\ref{eq:transLinearTS}) and for some $q > 0$ define $X_t^{(q)} = \bigoplus_{j=0}^{q}\psi_j \circ Z_{t-j}$ and $X_t^{(q)'} = \bigoplus_{j=q+1}^{\infty}\psi_j \circ Z_{t-j}$.
\citet[][Section 3.4]{mhatrecooley2020} say that $X_t^{(q)}$ converges to $X_t$ in tail ratio if $\lim \limits_{x \to \infty} \frac{\text{pr}
(X_t^{(q)'} > x)}{x^{-2}L(x)} = 0$ as $q\rightarrow \infty$.
The proofs below are extreme analogues to the steps in \citet[][Proposition 3.2.1]{brockwelldavis1991}, with tail ratio convergence replacing mean-square convergence.

\begin{lemma}
If $X_t$ is a tail stationary process in $\mathbb{V}$,
then
$$
    P_{\Bar{\text{sp}}\{X_j, t-n\le j\le t-1\}}X_{t}   \stackrel{\text{tail ratio}}{\rightarrow} P_{\Bar{\text{sp}}\{X_j, -\infty <j\le t-1\}}X_{t}, \; \; \text{as } n\to \infty. \nonumber
$$
\end{lemma}
\noindent{Proof}: 
Consider the transformed-linear combination 
$$
    \bigoplus^{\infty}_{j=n+1} b_j \circ X_{t-j}. \nonumber    
$$
As $X_t$ is a tail stationary process,
$$
    \text{pr} (X_t >x) \sim x^{-2}L(x)\sigma(0), \nonumber
$$
where $\sigma(0)=\sigma(X_t,X_t)$. 
As shown in \citet[Section 3.4]{mhatrecooley2020},
\begin{align} \label{eq:322}
    \text{pr} \left(\bigoplus^{\infty}_{j=n+1} b_j \circ X_{t-j} >x\right) & \sim x^{-2}L(x)\sigma(0)\sum^{\infty}_{j=n+1} |b_j|^2, \; \; \text{as }x\to \infty. 
\end{align}
Taking limit on both sides of (\ref{eq:322}) we get, as $x \to \infty$,
\begin{eqnarray}
    \lim \limits_{n \to \infty} \text{pr}\left(\bigoplus^{\infty}_{j=n+1} b_j \circ X_{t-j} >x\right) 
    & \sim&  \lim \limits_{n \to  \infty} x^{-2}L(x)\sigma(0)\sum^{\infty}_{j=n+1} |b_j|^2 = 0, \nonumber \\
   \implies \lim \limits_{x \to \infty} \frac{\text{pr}\left(\bigoplus^{\infty}_{j=n+1} b_j \circ X_{t-j} >x\right)}{x^{-2}L(x)} 
   & \to& 0, \; \; \text{as } n\to \infty. \nonumber
\end{eqnarray}
Thus
\begin{eqnarray}
    \bigoplus^n_{j=1} b_j \circ X_{t-j} & \stackrel{\text{tail ratio}}{\rightarrow} & \bigoplus^{\infty}_{j=1} b_j \circ X_{t-j}, \mbox{ as } n \rightarrow \infty, \mbox{ and } \nonumber \\
    P_{\Bar{\text{sp}}\{X_j, t-n\le j\le t-1\}}X_{t} &  \stackrel{\text{tail ratio}}{\rightarrow} & P_{\Bar{\text{sp}}\{X_j, -\infty <j\le t-1\}}X_{t}, \; \; \text{as } n\to \infty. \nonumber
\end{eqnarray}
\hspace*{0pt}\hfill $\square$

\noindent{Proof of Proposition \ref{prop:qcorrelated}}: For each $t$, let $\mathbb{M}_t$ be the closed transformed-linear subspace $\Bar{\text{sp}}\{X_s, s \le t\}$ of $\mathbb{V}$ and set
\begin{align} \label{eq:inn}
    Z_t  & = X_t \ominus P_{\mathbb{M}_{t-1}}X_t.
\end{align}
That is,
$$
    Z_t  = X_t \ominus \bigoplus^{\infty}_{j=1}b_j \circ X_{t-j}, \; \; a_j \in \mathbb{R}. \nonumber
$$
Thus, $Z_t \in \mathbb{M}_t$. By definition of $P_{\mathbb{M}_{t-1}}$, $P_{\mathbb{M}_{t-1}}X_t \in \mathbb{M}_{t-1}$ and $Z_t = X_t \ominus P_{\mathbb{M}_{t-1}}X_t \in \mathbb{M}_{t-1}^{\perp}$.
Thus $Z_s \in \mathbb{M}_s \subset \mathbb{M}_{t-1}$ and hence $\langle Z_s,Z_t \rangle = 0$ for $s<t$. Also, by Lemma 1
$$
    P_{\Bar{\text{sp}}\{X_s, t-n\le s\le t-1\}}X_{t}  \stackrel{\text{tail ratio}}{\rightarrow} P_{\mathbb{M}_{t-1}}X_{t}, \; \; \text{as } n\to \infty. \nonumber
$$
By stationarity and continuity of norm,
\begin{eqnarray*}
    \|Z_{t+1}\| & =& \|X_{t+1} \ominus P_{\mathbb{M}_t}X_{t+1}\| \nonumber \\
    & =& \lim \limits_{n \to \infty} \|X_{t+1} \ominus P_{\Bar{\text{sp}}\{X_s,s=t+1-n, \cdots, t\}}X_{t+1}\| \nonumber \\
    & =& \lim \limits_{n \to \infty} \|X_t \ominus P_{\Bar{\text{sp}}\{X_s,s=t-n, \cdots, t-1\}}X_t\| \nonumber \\
    & =& \|X_t \ominus P_{\mathbb{M}_{t-1}}X_t\| = \|Z_t\|. \nonumber
\end{eqnarray*}
Letting $c^2=\|Z_t\|^2$, $\{Z_t\}$ is a sequence of independent and tail stationary regularly-varying random variables with scale $c$, that is, $\text{Pr}(Z_t > x)/\{x^{-2}L(x)\} = c^2$.

\noindent{By} (\ref{eq:inn}),
$$
    X_{t-1} = Z_{t-1} \oplus P_{\mathbb{M}_{t-2}}X_{t-1}. \nonumber
$$
Consequently,
\begin{eqnarray}
    \mathbb{M}_{t-1} 
    & =& \Bar{\text{sp}}\{X_s, s\le t-1\} \nonumber \\
    & =& \Bar{\text{sp}}\{X_s, s< t-1, Z_{t-1}\} \nonumber  \\
    & =& \Bar{\text{sp}}\{X_s, s< t-q, Z_{t-q}, \cdots, Z_{t-1}\}. \nonumber
\end{eqnarray}
Therefore, $\mathbb{M}_{t-1}$ can be decomposed into two orthogonal subspaces, $\mathbb{M}_{t-q-1}$ and $\Bar{\text{sp}}\{Z_{t-q}, \cdots, Z_{t-1}\}$. Since $\gamma(h)=0$ for $|h|>q$, it follows that $X_t \perp \mathbb{M}_{t-q-1}$ and since $\Bar{\text{sp}}\{c^{-2}\circ Z_{t-q}, \cdots, c^{-2}\circ Z_{t-1}\}$ is an orthonormal set, by properties of projection mappings,
\begin{eqnarray} \label{prop3proof}
    P_{\mathbb{M}_{t-1}}X_t 
    & =& P_{\mathbb{M}_{t-q-1}}X_t \oplus P_{\Bar{\text{sp}}\{Z_{t-q}, \cdots Z_{t-1}\}}X_t \nonumber \\
    & =& 0 \; \oplus \; \left(c^{-2}\langle X_t,Z_{t-1} \rangle \right) \circ Z_{t-1} \;\oplus \cdots \oplus \; \left(c^{-2}\langle X_t,Z_{t-q} \rangle \right) \circ Z_{t-q} \nonumber \\
    & =& \theta_1 \circ Z_{t-1} \oplus \cdots \oplus \theta_q \circ Z_{t-q},
\end{eqnarray}
where $\theta_j:= c^{-2}\langle X_t,Z_{t-j} \rangle$, which by stationarity is independent of $t$ for $j=1, \cdots, q$. By (\ref{eq:inn}) and (\ref{prop3proof}),
$$
     X_t = Z_t \oplus \theta_1 \circ Z_{t-1} \oplus \cdots \oplus \theta_q \circ Z_{t-q}. \nonumber
$$
\hspace*{0pt}\hfill $\square$

\section{Modeling in Subset $\mathbb{V}_+$} \label{sec:subsetVplus}

We have defined $\mathbb{V}$ allowing for negative $\psi_j$'s as they are needed to have a vector space and to create orthogonal elements.
In this section we will show that the negative coefficients defining $\{X_t\}$ are academic in the sense that a time series that has negative coefficients is indistinguishable in terms of tail behavior from a time series that has zeroes in place of those coefficients.
Rather than being a problem, this will allow us to restrict our attention to a subset $\mathbb{V}_+$ for which the TPDF is equivalent to the inner product function. 

Consider a subset of $\mathbb{V}$ defined as
$$\mathbb{V}_+ = \{X_t : X_t = \bigoplus^{\infty}_{j=0}\psi_j \circ Z_{t-j}, \; \psi_j\ge 0, \sum^{\infty}_{j=0} \psi_j<\infty\}.$$
Proposition \ref{prop:vplus} below follows from the definition of TPDF.

\begin{proposition}\label{prop:vplus}
If a transformed-linear MA($\infty$) time series in $\mathbb{V}$ has TPDF $\sigma(h)$, then there exists a transformed-linear MA($\infty$) time series in subset $\mathbb{V}_+$ which has the same TPDF $\sigma(h)$, for all lag $h$.
\end{proposition}

\noindent{Proof}: Let $X_t = \bigoplus^{\infty}_{j=0}\psi_j \circ Z_{t-j} \in \mathbb{V}$.
The TPDF of $X_t$ is given by $\sigma(h) = \sum^{\infty}_{j=0} \psi^{(0)}_j \psi^{(0)}_{j+h}$, which is equal to the TPDF of $X_t^* = \bigoplus^{\infty}_{j=0}\psi_j^{(0)} \circ Z_{t-j} \in \mathbb{V}_+.$ 
\hspace*{0pt}\hfill $\square$

\noindent{In other words,} $X_t$ and $X_t^*$ are indistinguishable in terms of tail dependence.
$X_t$ and $X_t^*$ are also indistinguishable in terms of the tail ratio since recall that tail ratio is equal to $\sigma(0)$.
Furthermore, the TPDF gives full information for a stationary time series in $\mathbb{V}_+$ and unlike the inner product function, the TPDF is estimable.
Also, it can be clearly seen that $\gamma(h) = \sigma(h)$, for $X_t^*, X^*_{t+h} \in \mathbb{V}_+$, for all lag $h$.
Hence it is reasonable to restrict our attention to $\mathbb{V}_+$.

As the inner product is equivalent to the TPDF in $\mathbb{V}_+$, equation \ref{eq:bn} can be rewritten as 
\begin{align} \label{b_hat_v+}
    \hat{\boldsymbol{b}}_n = \Sigma_n^{-1}\boldsymbol{\sigma}_n.
\end{align}
where $\Sigma_n = \left[ \sigma(i-j) \right]^n_{i,j=1}$ and $\boldsymbol{\sigma}_n = \left[ \sigma(i) \right]^n_{i=1}$.

Also, if our time series is in $\mathbb{V}_+$, we can rewrite the equations (\ref{ia}) of the innovations algorithm in terms of the TPDF $\sigma(\cdot)$ instead of the inner product as
\begin{eqnarray} \label{vplusIA}
    \begin{cases}
      \nu_0 & = \sigma(0)\\
      \theta_{n,n-k} & = \nu_k^{-1}\left(\sigma(n-k) - \sum_{j=0}^{k-1}\theta_{k, k-j} \theta_{n, n-j} \nu_j \right), \quad k=0,1,...,n-1,\\
      \nu_n & = \sigma(0) - \sum_{j=0}^{n-1} \theta_{n,n-j}^2 \nu_j,
    \end{cases}.       
\end{eqnarray}
Rewriting Corollary \ref{prop:convergence} for $\mathbb{V}_+$, we get the following corollary.
\begin{corollary}\label{cor:convergence_vplus}
If $\{X_t\}$ is an invertible MA process in $\mathbb{V}_+$ with $\lim_{x \to \infty} \text{Pr}(Z_j > x)/\{x^{-2}L(x)\} = 1$, then as $n\to \infty$, \\
(i) $\nu_n \to 1$, \\
(ii) $\|X_n \ominus \hat{X}_n \ominus Z_n\|^2 \to 0$, and \\
(iii) $\theta_{nj} \to \psi_j, j=1,2,\cdots$; \; $\psi_j\ge 0$.
\end{corollary}

Also, rewriting Proposition \ref{prop:qcorrelated} for $\mathbb{V}_+$,
\begin{corollary}\label{cor:qcorrelated_vplus}
If $\{X_t\}$ is a regularly-varying tail stationary process in $\mathbb{V}_+$ with tail pairwise dependence function $\sigma(.)$ such that $\sigma(h)=0$ for $|h|>q$ and $\sigma(q)\ne 0$, then $\{X_t\}$ is an transformed-linear regularly-varying MA($q$) process, i.e. there exists a regularly-varying noise sequence $\{Z_t\}$ of independent and tail stationary $Z_t$'s such that
$$
    X_t = Z_t \oplus \theta_1 \circ Z_{t-1} \oplus \cdots \oplus \theta_q \circ Z_{t-q}. \nonumber
$$
\end{corollary}

The relation between the TPDF and the inner product gives an important result described in the following remark:
\begin{remark}
If $X_t$ is an MA($\infty$) time series in $\mathbb{V}$ and $X_t \notin \mathbb{V}_+$, then by Proposition \ref{prop:vplus} there exists $X_t^* \in \mathbb{V}_+$, obtained by applying the zero-operator on the coefficients of $X_t$, which has the same TPDF as $X_t$. Thus, the innovations algorithm applied to the TPDF of $X_t \in \mathbb{V}$ will give us the one-step predictors for the corresponding $X_t^* \in \mathbb{V}_+$.
\end{remark}

\section{Flexibility of the MA($\infty$) Class for Modeling}
\label{denseness}

In this section, we show that the class of MA($\infty$) models is a rich class for modeling.

\subsection{Richness of the MA($\infty$) Class in terms of the TPDF}


Given a valid TPDF (that is, a completely positive function) that converges to $0$ as lag increases, we can run the innovations algorithm to get the $\theta_{nj}$'s and $\nu_n$ as defined in (\ref{ia}).
If we apply the TPDF formula to these $\theta_{nj}$'s we will get a TPDF that gets arbitrarily close to the given TPDF.
In other words, if we consider random noise terms $\{Z_j\}$ that are Fr\`echet with $\alpha=2$ and scale $\sqrt{\nu_n}$, and generate a process applying the coefficients $\theta_{nj}$ to the $Z$'s, the TPDF of this generated process will be arbitrarily close to the given TPDF.
Thus, given any valid TPDF, we can run the transformed-linear innovations algorithm long enough to find a transformed-linear regularly-varying MA($\infty$) time series whose TPDF will get arbitrarily close to the given TPDF.
As such, the class of MA($\infty$) time series is rich in the class of possible TPDFs that converge to $0$.

To show this, first, we prove the result for a $q$-tail-dependent TPDF in the following corollary.

\begin{corollary}
If $\{X_t\}$ is any regularly-varying tail stationary process with TPDF $\sigma(.)$ such that $\sigma(h)=0$ for $|h|>q$ and $\sigma(q)\ne 0$, then as $n \to \infty$, the $\theta_{nj}$s generated from the transformed-linear innovations algorithm approach $\theta_1, \cdots, \theta_q$ of an MA($q$) whose TPDF matches the given TPDF.
\end{corollary}

\noindent{Proof:} Let us consider the form for the $\theta_{nj}$s given by the transformed-linear innovations algorithm in \ref{vplusIA}:
\begin{eqnarray}\label{theta1}
\theta_{n,n-k} = \nu_k^{-1}\left(\sigma(n-k) - \sum_{l=0}^{k-1}\theta_{k, k-l} \theta_{n, n-l} \nu_l \right), \quad k=0,1,...,n-1.
\end{eqnarray}
Rewriting (\ref{theta1}) by letting $h = n-k$,
\begin{eqnarray}\label{theta2}
\theta_{n,h} 
& =& \nu_{n-h}^{-1}\left(\sigma(h) - \sum_{l=0}^{n-h-1}\theta_{n-h, n-h-l} \theta_{n, n-l} \nu_l \right),  \nonumber \\
& =& \nu_{n-h}^{-1}\left(\sigma(h) - \sum_{l=n-h-q}^{n-h-1}\theta_{n-h, n-h-l} \theta_{n, n-l} \nu_l \right), \quad h=0,1,...,q,
\end{eqnarray}
since $\theta_{n,n-l} = 0$ for all $l = 0,1, \cdots, n-h-q-1$.

\noindent{Rewriting} (\ref{theta2}) by letting $ j= n-h-l$,
\begin{align}
  \label{theta3}
  \theta_{n,h} = \nu_{n-h}^{-1}\left(\sigma(h) - \sum_{j=1}^{q}\theta_{n-h, j} \theta_{n, j+h} \nu_{n-h-j} \right).
\end{align}
As $n \to \infty$, let $\theta_{n,h} \to \theta_h$ and $\nu_n \to c^2$. Then (\ref{theta3}) becomes,
\begin{align}
  \label{theta4}
  \theta_{h} &= c^{-2}\left(\sigma(h) - \sum_{j=1}^{q}\theta_j \theta_{j+h} c^{2} \right).
\end{align}
Rearranging (\ref{theta4}),
\begin{eqnarray}\label{theta5}
\sigma(h) 
&=& \theta_{h} c^{2} + \sum_{j=1}^{q}\theta_j \theta_{j+h} c^{2} \nonumber \\
&=& c^{2}\sum_{j=0}^{q}\theta_j \theta_{j+h}, \quad h=0,1,...,q,
\end{eqnarray}
which is the form for the TPDF at lag $h$ for a regularly-varying tail stationary MA($q$) with $\theta_j \ge 0$ for $j=0,1, \cdots, q$ and $\lim_{x \to \infty} \text{pr}(Z_j > x)/\{x^{-2}L(x)\} = c^2$. Thus, the TPDF of this MA($q$) matches the given TPDF $\sigma(h)$. \hfill{$\square$}

We are extending the above result to the MA($\infty$) case.

\subsection{Transformed-Linear Wold Decomposition}

If the TPDF of a time series does not converge to $0$, analogous to the Wold decomposition discussed in \citet{brockwelldavis1991} we can decompose the time series into an MA($\infty$) process and a deterministic process.
Following \citet{brockwelldavis1991} and \citet{sargent1979} we prove our Transformed-Linear Wold Decomposition as follows:
\begin{theorem}[The Transformed-Linear Wold Decomposition]\label{thm:wold}
If $c^2 = \|X_{n+1} \ominus \hat{X}_{n+1}\|^2>0$, then $X_t$ can be expressed as
\begin{eqnarray} \label{eq:wold}
    X_t = \bigoplus_{j=0}^{\infty} \psi_j \circ Z_{t-j} \; \oplus \; U_t,
\end{eqnarray}
where
\begin{eqnarray}
    & \text{(i) } \psi_0=1 \text{ and } \sum_{j=0}^{\infty}\psi_j^2 < \infty, \nonumber \\
    & \text{(ii) } \{Z_t\} \text{ is a sequence of independent and tail stationary regularly-varying random variables with scale } c, \nonumber \\
    & \text{(iii) } Z_t \in \mathbb{V}^t \text{ for each } t\in \mathbb{Z},  \nonumber \\
    & \text{(iv) } \langle Z_t, U_s \rangle =0 \text{ for all } s,t \in \mathbb{Z}, \nonumber \\
    & \text{(v) } \{U_t\} \text{ is deterministic.} \nonumber
\end{eqnarray}
The sequences $\{\psi_j\}$, $\{Z_j\}$, and $\{U_j\}$ are uniquely determined by (\ref{eq:wold}) and the conditions (i) - (v).
\end{theorem}

\noindent {Proof:} Consider the sequences
\begin{align} \label{eq:Z_t}
    Z_t = X_t \ominus P_{\mathbb{V}^{t-1}}X_t, 
\end{align} 
\begin{align} \label{eq:psij}
    \psi_j = c^{-2}\langle X_t , Z_{t-j} \rangle,
\end{align}
\begin{eqnarray} \label{eq:u_t}
    U_t = X_t \ominus \bigoplus_{j=0}^{\infty}\psi_j \circ Z_{t-j}.
\end{eqnarray}

That is,
$$
    Z_t = X_t \ominus \bigoplus^{\infty}_{j=1}a_j \circ X_{t-j}, \; \; a_j \in \mathbb{R}, \; j = 1, \cdots, t-1. \nonumber
$$
Thus, $Z_t \in \mathbb{V}^t$, establishing (iii).
By definition of $P_{\mathbb{V}^{t-1}}$, $P_{\mathbb{V}^{t-1}}X_t \in \mathbb{V}^{t-1}$ and $Z_t = X_t \ominus P_{\mathbb{V}^{t-1}}X_t \in \mathbb{V}^{{t-1}{\perp}}$.
Thus,$$Z_t \in \mathbb{V}^{{t-1}{\perp}} \subset \mathbb{V}^{{t-2}{\perp}} \subset \cdots$$
Hence for $s<t$, $\langle Z_s,Z_t \rangle = 0$.
By Lemma 1
$$
    P_{\Bar{\text{sp}}\{X_s, t-n\le s\le t-1\}}X_{t}   \stackrel{\text{tail ratio}}{\rightarrow} P_{\mathbb{V}^{t-1}}X_{t}, \; \; \text{as } n\to \infty. \nonumber
$$
By stationarity and continuity of norm,
\begin{eqnarray}
    \|Z_{t+1}\| &=& \|X_{t+1} \ominus P_{\mathbb{V}^t}X_{t+1}\| \nonumber \\
    &=& \lim \limits_{n \to \infty} \|X_{t+1} \ominus P_{\Bar{\text{sp}}\{X_s,s=t+1-n, \cdots, t\}}X_{t+1}\| \nonumber \\
    &=& \lim \limits_{n \to \infty} \|X_t \ominus P_{\Bar{\text{sp}}\{X_s,s=t-n, \cdots, t-1\}}X_t\| \nonumber \\
    &=& \|X_t \ominus P_{\mathbb{V}^{t-1}}X_t\| = \|Z_t\|. \nonumber
\end{eqnarray}
Letting $c^2=\|Z_t\|^2$, $\{Z_t\}$ is a sequence of independent and tail stationary regularly-varying random variables with scale $c$, thus establishing (ii).

By equation (\ref{prop3proof}) in the proof of Proposition \ref{prop:qcorrelated},
$$
    P_{\Bar{\text{sp}}\{Z_j, j\le t\}}X_{t} = \sum_{j=0}^{\infty} \psi_j \circ Z_{t-j}, \nonumber
$$
where $\psi_j=c^{-2}\langle X_t , Z_{t-j} \rangle$ and $\sum_{j=0}^{\infty} \psi_j^2 < \infty$. By stationarity, the coefficients $\psi_j$ are independent of $t$.
Also,
$$
    \psi_0 = c^{-2}\langle X_t, X_t \ominus P_{\mathbb{V}^{t-1}}X_t \rangle = c^{-2} \| X_t \ominus P_{\mathbb{V}^{t-1}}X_t \|^2  = c^{-2} \|Z_t\|^2 = 1, \nonumber
$$
thus proving (i). From equation (\ref{eq:psij}) and (\ref{eq:u_t}), for $s \le t$,
\begin{eqnarray} 
    \langle U_t, Z_s \rangle 
    & =& \left \langle X_t \ominus \bigoplus_{j=0}^{\infty}\psi_j \circ Z_{t-j}, \; Z_s \right \rangle \nonumber \\
    & =& \left \langle X_t, Z_s \right \rangle - \left \langle \bigoplus_{j=0}^{\infty}\psi_j \circ Z_{t-j}, \; Z_s \right \rangle \nonumber \\
    & =& \left \langle X_t, Z_s \right \rangle - \psi_{t-s} \left \langle Z_s, Z_s \right \rangle \nonumber \\
    & =& \left \langle X_t, Z_s \right \rangle - \|Z_s\|^{-2}\left \langle X_t, Z_s \right \rangle \|Z_s\|^{2} \nonumber \\
    & =& 0. \nonumber
\end{eqnarray}
In addition, if $s>t$, $Z_s \in \mathbb{V}^{{s-1}{\perp}} \subset \mathbb{V}^{{t}{\perp}}$. But $U_t \in \mathbb{V}^t$. Hence $\langle U_t, Z_s \rangle = 0$ for $s>t$. Thus (iv) is proved.

Since $U_t$ is orthogonal to $Z_t$, $U_t \in \mathbb{V}^{t-1}$, that is $U_t$ can be predicted perfectly from previous $X$'s.
To see this clearly, consider the projection of $U_t$ on $\mathbb{V}^{t-1}$ to get
$$
P_{\mathbb{V}^{t-1}}U_t = P_{\mathbb{V}^{t-1}}X_t \ominus P_{\mathbb{V}^{t-1}}\bigoplus_{j=0}^{\infty}\psi_j \circ Z_{t-j} \nonumber \\
= P_{\mathbb{V}^{t-1}}X_t \ominus \; \bigoplus_{j=1}^{\infty}\psi_j \circ Z_{t-j}, \nonumber
$$
since $P_{\mathbb{V}^{t-1}}Z_t = 0$ and $P_{\mathbb{V}^{t-1}}Z_{t-k} = Z_{t-k}$ for $k \ge 1$.
Transformed-linearly subtracting above equation from (\ref{eq:u_t}) gives
$$U_t \ominus P_{\mathbb{V}^{t-1}}U_t = \left( X_t \ominus P^{\mathbb{V}^{t-1}}X_t \right) \; \ominus \; \psi_0\circ Z_t = 0_{\mathbb{V}},$$
since the one-step ahead prediction error for $X_t$ is $\psi_0\circ Z_t$.
Hence, $U_t=P_{\mathbb{V}^{t-1}}U_t$.
In general,
$$
P_{\mathbb{V}^{t-k}}U_t = P_{\mathbb{V}^{t-k}}X_t \ominus \; \bigoplus_{j=k}^{\infty}\psi_j \circ Z_{t-j}. \nonumber
$$
Transformed-linearly subtracting above equation from (\ref{eq:u_t}) gives
$$U_t \ominus P_{\mathbb{V}^{t-k}}U_t = \left( X_t \ominus P_{\mathbb{V}^{t-k}}X_t \right) \; \ominus \; \bigoplus_{j=0}^{k-1}\psi_j \circ Z_{t-j} = 0_{\mathbb{V}},$$
since the k-step ahead prediction error for $X_t$ is $\bigoplus_{j=0}^{k-1}\psi_j \circ Z_{t-j}$.
Thus $\{U_t\}$ is deterministic as it can be predicted from past $X$'s.\hfill $\square$

\subsection{Simulation Study}
\label{simulation_study}

\begin{figure}
\includegraphics[width=3.2 in]{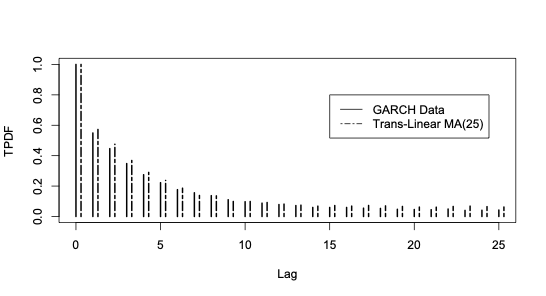}
\includegraphics[width=3.2 in]{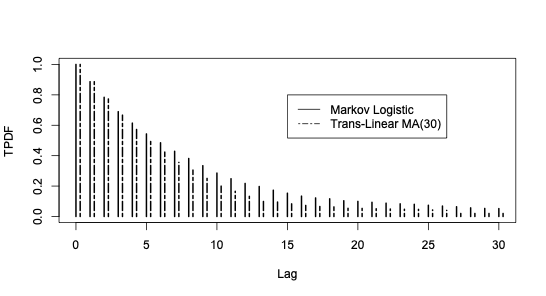}
  \caption{Left panel:  estimated TPDF's for the original GARCH data and for the data simulated from the fitted transformed-linear MA(25) model.  Right panel:  same but for the original Markov logistic model and the fitted transformed-linear MA(30).}
  \label{fig:simStudyTPDF}
\end{figure}

We conduct a simulation study that corroborates the richness of the class of transformed-linear regularly-varying MA($\infty$) models.
We simulate data from two models.
The first model is a GARCH(1,1) process (\citet{bollerslev1986}) with Gaussian noise terms and parameters $\alpha_0=0.2$, $\alpha_1=0.5$, and $\beta_1=0.3$.
We consider the time series of absolute values of this GARCH process, which we denote $x_t^{(orig)}$
A chi-plot (not shown) for the upper tail shows asymptotic dependence with $\hat{\chi}(1)\approx 0.34$.
The Hill estimator (\citet{hill1975}) at the empirical $0.99$ quantile of this transformed data gives an estimate $\hat{\alpha}=3.27$ of the tail index.
The scale is estimated to be $\hat{c}=0.47$.
We further transform the data into $x_t=\hat{c}^{-1/2}(x_t^{(orig)})^{\hat{\alpha}/2}$ so that our marginal now can be assumed to have. $\alpha=2$ and $\sigma(0)=1$.
As discussed in \citet{mhatrecooley2020}, preprocessing the data to have $\sigma(0)=1$ allows us to reduce bias in TPDF estimation.
Note that by doing this the noise terms $Z_j$ are no longer such that $\sigma_{Z_j}(0) = 1$.
As done in \citet{mhatrecooley2020}, to reduce bias in TPDF estimation, we subtract off the mean of the transformed data and replace the negative observations with $0$.
We estimate the TPDF up to 500 lags using data whose radial components exceeds the $0.99$ quantile.
The squared distance of prediction $\nu_n$ converges to 0.65.
Running the innovations algorithm on the estimated TPDF gives us converged $\theta$ estimates of an MA model.
We consider the first 25 $\hat{\theta}$'s as we deem the $\theta$ estimates to be negligible beyond that.
We then generate Fr\'echet noise terms with $\alpha=2$ and scale $\sqrt{\nu_n} = \sqrt{0.65}$, and simulate a transformed-linear regularly-varying MA($25$) time series using the estimated $\theta$'s from the innovations algorithm.
We then back-transform the simulated time series to the original marginals.
The average difference for the first 25 lags between the estimated TPDF from the original GARCH data and the estimated TPDF from our fitted model is $-0.01$ (se = $0.01$).
Figure \ref{fig:simStudyTPDF} shows the estimated TPDF's for both the GARCH model data and the data simulated from the fitted MA(25).

The second model is a first-order Markov chain such that each pair of consecutive observations has a bivariate logistic distribution with a dependence parameter of 0.4 and common unit-Fr\'echet marginals (refer \citet{smithetal}).
A chi-plot (not shown) for the upper tail shows asymptotic dependence with $\hat{\chi}(1)\approx 0.7$.
We perform the square-root transformation on the data so that $\alpha=2$.
Following the same process as for the first model, we simulate a transformed-linear regularly-varying MA($30$) time series using the estimated $\theta$'s from the innovations algorithm and back-transform the simulated time series to the original marginals.
The average difference for the first 30 lags between the estimated TPDF from the original logistic data and the estimated TPDF from our fitted model is $0.05$ (se = $0.03$).

Table \ref{table_run_dense} gives the average length of run above higher quantiles for the original time series data and the fitted time series data (using coefficient estimates from the innovations algorithm) for the GARCH model and the logistic model. Table \ref{table_sum_dense} gives the higher quantiles for sum of three consecutive time series terms. The fitted models seem to produce reasonable estimates of these tail quantities. Interestingly, for the fitted MA model in the logistic case (last column in Table \ref{table_run_dense}), there is an increasing trend in the estimates and we do not see the ``threshold stability'' as exhibited by the fitted GARCH model.
Tables \ref{table_run_dense} and \ref{table_sum_dense} show that a model which captures the level of dependence in the TPDF can adequately estimate quantities of interest like length of run or sums of consecutive terms, despite not fully characterizing the angular measure of the true model.
In fact, lag-1 scatterplots of the target and fitted models (available in \cite{mhatre2022}) show distinct differences in the nature of the bivariate dependence structure.

\begin{table}[h]
\fontsize{10}{10}\selectfont
 \begin{center}
    \caption{Average length (standard error) of run above a threshold for the simulation study}
    \label{table_run_dense}
        \begin{tabular}{lccccc}
        \hline
        Threshold & \multicolumn{2}{c}{GARCH} &  &\multicolumn{2}{c}{Logistic} \\ 
                             quantile   & \multicolumn{1}{c}{Original} & \multicolumn{1}{c}{Fitted}   & \multicolumn{1}{c}{} & \multicolumn{1}{c}{Original}  & \multicolumn{1}{c}{Fitted} \\ \hline
        0.95 & \multicolumn{1}{c}{1.57 (0.02)} & \multicolumn{1}{c}{1.71 (0.03)}  & \multicolumn{1}{c}{} & \multicolumn{1}{c}{3.02 (0.09)} & \multicolumn{1}{c}{3.88 (0.11)} \\
        0.98 & \multicolumn{1}{c}{1.57 (0.03)} & \multicolumn{1}{c}{1.60 (0.04)}  & \multicolumn{1}{c}{} & \multicolumn{1}{c}{3.21 (0.15)} & \multicolumn{1}{c}{3.87 (0.16)}\\ 
        0.99 & \multicolumn{1}{c}{1.56 (0.05)} & \multicolumn{1}{c}{1.62 (0.06)}  & \multicolumn{1}{c}{} & \multicolumn{1}{c}{3.44 (0.22)} & \multicolumn{1}{c}{4.46 (0.26)}\\ 
        0.995 & \multicolumn{1}{c}{1.56 (0.06)} & \multicolumn{1}{c}{1.66 (0.09)}  & \multicolumn{1}{c}{} & \multicolumn{1}{c}{3.45 (0.29)} & \multicolumn{1}{c}{4.63 (0.35)}\\ 
        0.999 & \multicolumn{1}{c}{1.47 (0.11)} & \multicolumn{1}{c}{1.52 (0.15)}  & \multicolumn{1}{c}{} & \multicolumn{1}{c}{2.86 (0.48)} & \multicolumn{1}{c}{4.76 (0.76)}\\ \hline
        \end{tabular}
 \end{center}
\end{table}
\begin{table}[h]
 \fontsize{10}{10}\selectfont
 \begin{center}
    \caption{Quantiles for sum (standard error) of twelve consecutive terms for the simulation study}
    \label{table_sum_dense}
        \begin{tabular}{lccccc}
        \hline
        Quantile & \multicolumn{2}{c}{GARCH} &  &\multicolumn{2}{c}{Logistic} \\ 
                              & \multicolumn{1}{c}{Original} & \multicolumn{1}{c}{Fitted}   & \multicolumn{1}{c}{} & \multicolumn{1}{c}{Original}  & \multicolumn{1}{c}{Fitted} \\ \hline
0.95&	16.23	 (0.20) &	17.49	 (0.16) & & 	223.19	 (8.29) &	233.21	 (8.26) \\
0.98&	21.36	 (0.31) &	22.05	 (0.34) & & 	565.32	 (37.79) &	589.23	 (37.77) \\
0.99&	25.72	 (0.68) &	26.72	 (0.58) & & 	1252.96	 (184.83) &	1250.13	 (118.99) \\
0.995&	32.08	 (1.24) &	32.09	 (1.03) & & 	2789.53	 (420.92) &	2672.24	 (338.49) \\
0.999&	50.67	 (3.67) &	48.90	 (2.98) & & 	12833.18	 (3736.45) &	12732.06	 (3705.63) \\
 \hline
        \end{tabular}
 \end{center}
\end{table}



\section{Prediction Error}
\label{sec:pred_error}

We now return our attention to prediction and investigate the problem of assessing prediction uncertainty.
Because the geometry of regular variation is very different from the elliptical geometry assumed in many non-extreme settings, we need to deal with uncertainty in prediction differently.

\subsection{Completely Positive Decomposition of the Prediction TPDM}
\label{sec:CPD}

The squared distance of prediction, $\nu_n$, is the analogue to mean square prediction error.
In the finite-dimensional multivariate case \citet{leecooley2021} have shown that squared-norm prediction error $\nu_n$ is not useful to construct a prediction interval in the polar geometry of regular variation because the magnitude of error is dependent on the magnitude of the predicted value. We follow \citet{leecooley2021} and apply their method to construct prediction intervals when $\hat{X}_{n+1}$ is large.

The tail dependence between $\hat{X}_{n+1}$ and $X_{n+1}$ can be characterized by the bivariate angular measure $H_{\hat{X}_{n+1}, X_{n+1}}$. 
As shown in \citet{leecooley2021}, the dependence of $H_{\hat{X}_{n+1}, X_{n+1}}$ is summarized by the prediction TPDM
\begin{eqnarray} \label{sigma_pred}
\Sigma_{\hat{X}_{n+1}, X_{n+1}} = \begin{bmatrix} \boldsymbol{\sigma}_n^T\Sigma_n^{-1}\boldsymbol{\sigma}_n & \boldsymbol{\sigma}_n^T\Sigma_n^{-1}\boldsymbol{\sigma}_n \\
\boldsymbol{\sigma}_n^T\Sigma_n^{-1}\boldsymbol{\sigma}_n & \sigma(0) \end{bmatrix},
\end{eqnarray}
where $\Sigma_n = \left[ \sigma(i-j) \right]^n_{i,j=1}$ and $\boldsymbol{\sigma}_n = \left[ \sigma(i) \right]^n_{i=1}$.
Since $\Sigma_{\hat{X}_{n+1}, X_{n+1}}$ is a $2\times 2$ completely positive matrix, given any $q_*\ge 2$, there exist nonnegative matrices $B \in \mathbb{R}^{2\times q_*}$ such that $BB^T = \Sigma_{\hat{X}_{n+1}, X_{n+1}}$.
For feasible computation, \citet{leecooley2021} choose a moderate $q_*$ and apply the algorithm in \citet{groetznerdur2020} repeatedly to get $n_{decomp}$ nonnegative $B^{(k)}$ matrices, $k=1, \cdots, n_{decomp}$, such that $B^{(k)}B^{(k)^T}=\Sigma_{\hat{X}_{n+1}, X_{n+1}}$ for all $k$.
Then,
$$\hat{H}_{\hat{X}_{n+1}, X_{n+1}} = n_{decomp}^{-1}\sum_{k=1}^{n_{decomp}} \sum_{j=1}^{q_*}\|b^{(0)}_{kj}\|_2^2 \; \delta_{b^{(0)}_{kj}/\|b^{(0)}_{kj}\|_2}(\cdot),$$ where $b_{kj}$ is the $j^{\text{th}}$ column of the $k^{\text{th}}$ matrix $B^{(k)}$ and $\delta$ is the Dirac mass function, and
$$\Sigma_{\hat{H}}=n_{decomp}^{-1}\sum_{k=1}^{n_{decomp}}B^{(k)}B^{(k)^T} = \Sigma_{\hat{X}_{n+1}, X_{n+1}}.$$
Defined this way, $\hat{H}_{\hat{X}_{n+1}, X_{n+1}}$ consists of $n_{decomp}\times q_*$ point masses.

As in Section \ref{simulation_study}, we simulate 100,000 random observations of a first order Markov chain model such that each pair of consecutive observations has a bivariate logistic distribution with dependence parameter of 0.4 and common unit-Fr\'echet marginals.
We perform the square-root transformation on the data so that $\alpha=2$.
We consider the first 70,000 observations as training data and the remaining as test data.
In Section \ref{simulation_study} we fitted a transformed-linear regularly-varying MA($30$) time series to data simulated from the logistic model because the converged innovations algorithm gave negligible $\theta$ estimates beyond $\theta_{30}$.
Hence we consider the problem of predicting any observation $X_{n+1}$, $n \ge 30$, based on the previous $30$ observations.
Let us denote this predicted observation as $\hat{X}_{n+1|n:n-29}$.
Using equation (\ref{b_hat_v+}) we obtain $\hat{\boldsymbol{b}}$ and the prediction TPDM $\Sigma_{\hat{X}_{n+1|n:n-29}, X_{n+1}}$ is obtained from equation (\ref{sigma_pred}).
We apply the algorithm given in \citet{groetznerdur2020} repeatedly to compute $2 \times 5$ matrices $B^{(k)}$, $k=1, \cdots, 100$, each of which is a completely positive decomposition of $\Sigma_{\hat{X}_{n+1|n:n-29}, X_{n+1}}$.
Thus our estimated angular measure $\hat{H}_{\hat{X}_{n+1|n:n-29}, X_{n+1}}$ has $500$ point masses.
The $0.025$ and $0.975$ quantiles of $\hat{H}_{\hat{X}_{n+1|n:n-29}, X_{n+1}}$ give us a $95\%$ joint region.
The left panel of Figure \ref{joint_region} gives $95\%$ joint region on the 30,000 test data.
Thresholding at the $0.95$ quantile of $\|\hat{X}_{n+1|n:n-29}, X_{n+1}\|$, $99.6\%$ of the large data points fall within this joint region.

\subsection{Conditional Prediction Intervals}
\label{sec:prediction_intervals}

\begin{figure}[h]
\centerline{\includegraphics[width=0.5\textwidth, trim= 0cm 0cm 0cm 1.5cm, clip]{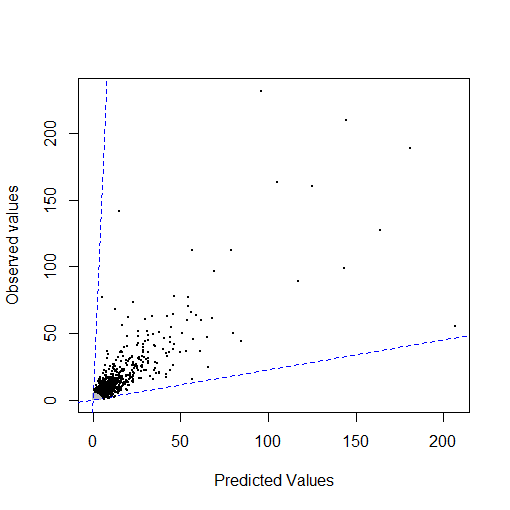}
    \includegraphics[width=0.5\textwidth, trim= 0cm 0cm 0cm 1.5cm, clip]{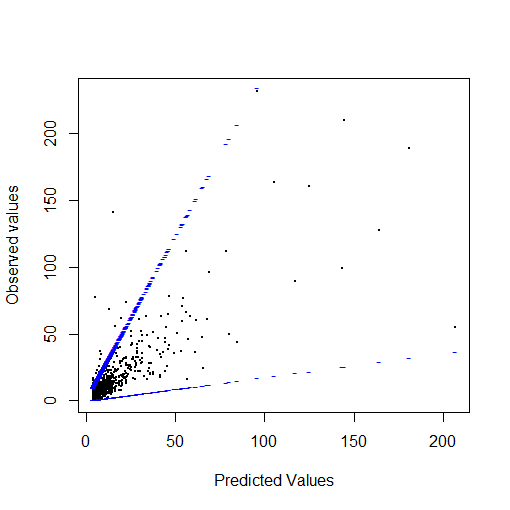}}
     \caption{Scatterplot of the Logistic model test data, on the transformed Fr\'echet scale with $\alpha=2$, with the estimated 95\% joint prediction region  (left panel). 95\% conditional prediction intervals given each large value of $\hat{X}_{n+1|n:n-29}$ of the Logistic model test data, on the transformed Fr\'echet scale with $\alpha=2$ (right panel).}
  \label{joint_region}
\end{figure}

The conditional density of $X_2|X_1=x_1$ if $x_1$ is large is given in \citet{leecooley2021} as approximately
\begin{align} \label{conditional}
    f_{X_2|X_1}(x_2|x_1) = 2c^{-1} \|(x_1,x_2)\|^{-5}_2 x_2 h\left( \frac{(x_1,x_2)}{\|(x_1,x_2)\|_2} \right),
\end{align}
where $c=\int^{\infty}_0 2\|(x_1,x_2)\|^{-5}_2 x_2 h \left( \frac{(x_1,x_2)}{\|(x_1,x_2)\|_2} \right)\text{d}x_2$.
We obtain an estimate of the conditional density of $X_{n+1}$ given a large value of $\hat{X}_{n+1}$ using equation (\ref{conditional}).
The angular density $h$ is estimated through a kernel density estimate of $\hat{H}_{\hat{X}_{n+1}, X_{n+1}}$.
The $0.025$ and $0.975$ quantiles of the estimated conditional density in equation (\ref{conditional}) give us a $95\%$ conditional prediction interval.
The right panel of Figure \ref{joint_region} gives the scatterplot after thresholding the test data at the $0.95$ quantile of $\hat{X}_{n+1|n:n-29}$, along with the $95\%$ conditional prediction intervals.
These prediction intervals have a coverage rate of $0.975$.

\section{Application to Santa Ana Winds}
\label{sec:app_prediction}

We return to the March AFB hourly windspeed data that we fitted a transformed-linear regularly-varying ARMA(1,1) time series to in \citet{mhatrecooley2020}.
We marginally transformed the data to have regularly-varying tails with $\alpha = 2$ and $\sigma(0) = 1$.
After subtracting off the mean of the transformed data and replacing the negative observations with $0$, we estimate the TPDF up to 500 lags using data whose radial components exceeds the $0.99$ quantile.
Running the innovations algorithm on the estimated TPDF gives us converged $\theta$ estimates of an MA model.
The squared distance of prediction, $\nu_n$, converges to 
We consider the first 40 $\hat{\theta}$'s since the $\theta$ estimates are negligible beyond that.
We then generate Fr\'echet noise terms, with $\alpha=2$ and scale $\sqrt{\nu_n} = \sqrt{0.65}$, and simulate a transformed-linear regularly-varying MA($40$) time series using the estimated $\theta$'s from the innovations algorithm and back transform the simulated time series to the original marginals.
The average difference between the estimated TPDF from the original windspeed anomalies data and the estimated TPDF from our fitted model is $-0.02$ (se = $0.01$).
Table \ref{table_run_dense_ws} gives average length of run above higher quantiles for the original windspeed anomalies time series data and the fitted time series data (using coefficient estimates from the innovations algorithm). Table \ref{table_sum_dense_ws} gives the higher quantiles for sum of three consecutive time series terms. The fitted models seem to produce reasonable estimates of these tail quantities.
\begin{table}[h]
\fontsize{10}{10}\selectfont
 \begin{center}
    \caption{Average length (standard error) of run above a threshold for the windspeed data}
    \label{table_run_dense_ws}
        \begin{tabular}{lcc}
        \hline
        Threshold & Original & Fitted \\ 
                             quantile   &  &\\ \hline
        0.95 & \multicolumn{1}{c}{2.43 (0.06)} & \multicolumn{1}{c}{2.32 (0.07)} \\
        0.98 & \multicolumn{1}{c}{2.35 (0.09)} & \multicolumn{1}{c}{2.27 (0.11)}\\ 
        0.99 & \multicolumn{1}{c}{2.10 (0.10)} & \multicolumn{1}{c}{2.46 (0.18)}\\ 
        0.995 & \multicolumn{1}{c}{1.77 (0.11)} & \multicolumn{1}{c}{2.35 (0.23)}\\ 
        0.999 & \multicolumn{1}{c}{1.40 (0.10)} & \multicolumn{1}{c}{2.04 (0.33)}\\ \hline
        \end{tabular}
 \end{center}
\end{table}
\begin{table}[h]
\fontsize{10}{10}\selectfont
 \begin{center}
    \caption{Quantiles for sum (standard error) of twelve consecutive terms for the windspeed data}
    \label{table_sum_dense_ws}
        \begin{tabular}{lcc}
        \hline
        Quantile & Original & Fitted \\ \hline
0.95&	27.70	 (0.72) &	28.44	 (0.52) \\
0.98&	43.74	 (1.19) &	42.93	 (1.07) \\
0.99&	56.65	 (1.66) &	54.82	 (1.62) \\
0.995&	69.67	 (2.11) &	68.91	 (2.50) \\
0.999&	91.89	 (3.60) &	97.93	 (4.49) \\
  \hline
        \end{tabular}
 \end{center}
\end{table}

\begin{figure}[h]
    \centerline{\includegraphics[width=0.33\textwidth, trim= 0cm 0cm 0cm 1.5cm, clip]{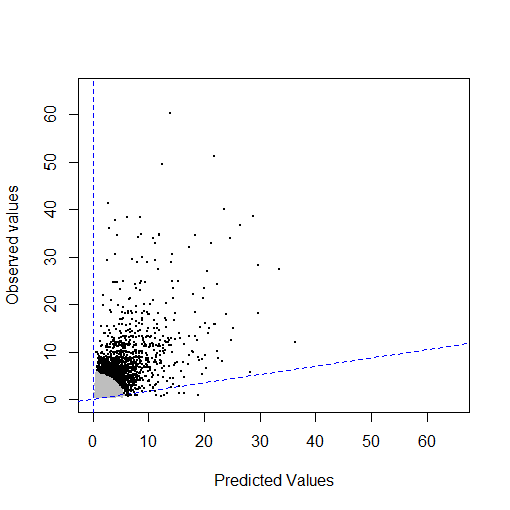}
    \includegraphics[width=0.33\textwidth, trim= 0cm 0cm 0cm 1.5cm, clip]{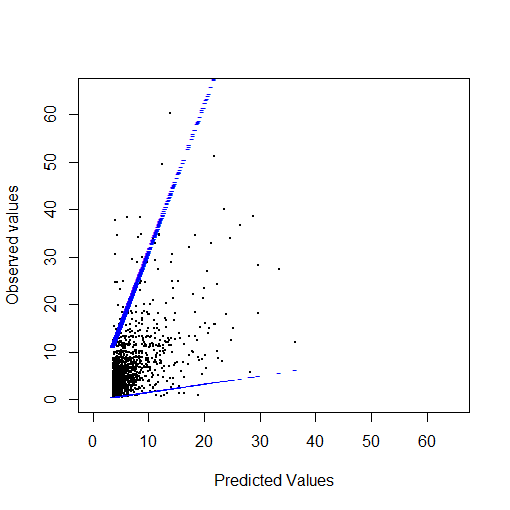}
    \includegraphics[width=0.33\textwidth, trim= 0cm 0cm 0cm 1.5cm, clip]{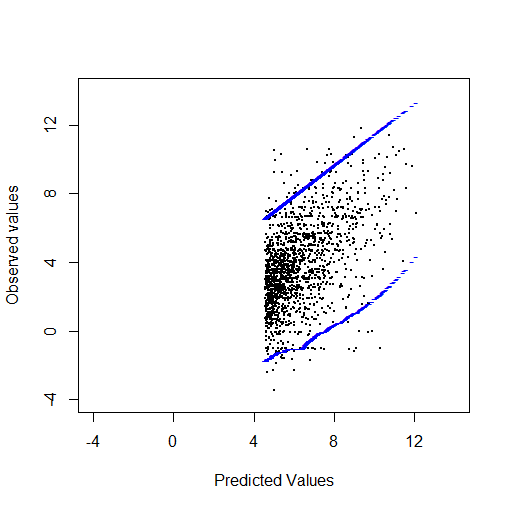}}
  \caption{Scatterplot of the windspeed anomalies test data on the Fr\'echet scale with the estimated 95\% joint prediction region (left panel). 95\% conditional prediction intervals given each large value of $\hat{X}_{n+1|n:n-39}$ of the windspeed anomalies test data on the Fr\'echet scale (center panel). 95\% conditional prediction intervals given each large value of $\hat{X}_{n+1|n:n-39}$ of the windspeed anomalies test data on the original scale (right panel).}
  \label{joint_region_ws}
\end{figure}

We now perform prediction of the windspeed anomalies time series.
Out of the 103,630 observations, we consider the first 70,000 observations as training data and the remaining as test data.
We consider the problem of predicting an observation $X_{n+1}$, $n \ge 40$, based on the previous 40 observations.
Let us denote this predicted observation as $\hat{X}_{n+1|n:n-39}$.
Using equation (\ref{b_hat_v+}) we obtain $\hat{\boldsymbol{b}}$ and the prediction TPDM $\Sigma_{\hat{X}_{n+1|n:n-39}, X_{n+1}}$ is obtained from equation (\ref{sigma_pred}).
We apply the algorithm given in \citet{groetznerdur2020} repeatedly to compute $2 \times 5$ matrices $B^{(k)}$, $k=1, \cdots, 100$, each of which is a completely positive decomposition of $\Sigma_{\hat{X}_{n+1|n:n-39}, X_{n+1}}$.
Thus our estimated angular measure $\hat{H}_{\hat{X}_{n+1|n:n-39}, X_{n+1}}$ has $500$ point masses.
The $0.025$ and $0.975$ quantiles of $\hat{H}_{\hat{X}_{n+1|n:n-39}, X_{n+1}}$ give us a $95\%$ joint region.
The left panel of Figure \ref{joint_region_ws} gives $95\%$ joint region on the test data.
Thresholding at the $0.95$ quantile of $\|\hat{X}_{n+1|n:n-39}, X_{n+1}\|$, $98.99\%$ of the large data points fall within this joint region.

We obtain an estimate of the conditional density of $X_{n+1}$ given a large value of $\hat{X}_{n+1}$ using equation (\ref{conditional}).
The $0.025$ and $0.975$ quantiles of the estimated conditional density in equation (\ref{conditional}) give us a $95\%$ conditional prediction interval.
The center panel of Figure \ref{joint_region_ws} gives the scatterplot after thresholding at the $0.95$ quantile of $\hat{X}_{n+1|n:n-39}$ of the test data along with the $95\%$ conditional prediction bounds.
These prediction intervals have a coverage rate of $0.96$.
The right panel of Figure \ref{joint_region_ws} gives the prediction intervals on the original scale of the anomalies obtained by taking the inverse of the marginal transformation.
We compare our prediction intervals to the standard Gaussian method.
We transform the marginal of the original windspeed anomalies data to be standard normal and estimate the ACVF.
We use the estimated covariance matrix to find the best linear unbiased predictor and to estimate the MSPE.
We then create 95\% Gaussian prediction intervals from the estimated MSPE and get a coverage rate of 0.94.
For the windspeed anomalies data, our prediction intervals do not show significant advantage over the standard Gaussian based prediction intervals because our data is not too heavy-tailed, resulting into a negligible difference between 
the corresponding predicted weight vectors $\hat{\boldsymbol{b}}$.
We are investigating a heavy-tailed precipitation data set where we suspect the difference between the 
corresponding predicted weight vectors will be more significant.


\section{Summary}

This paper extends the work of \cite{mhatrecooley2020}, which introduced the transformed linear time series models.
That paper laid out the foundation of these models, introducing the notion of tail stationarity which enables characterization of tail dependence in time series by the TPDF.
\cite{mhatrecooley2020} introduced the transformed linear backshift operator which allows the AR and ARMA models to be defined as transformed linear solutions to the model equations.
They showed via example that the transformed linear models were more flexible than nonnegative traditional linear models in capturing tail dependence.


This paper introduces the transformed linear innovations algorithm with the aim of performing prediction when the previous time series terms are large.
Not only does the innovations algorithm enable iterative prediction, but it
also provides us a tool demonstrating the richness of the class of models.
Perhaps the most important result is that one can fit a transformed-linear regularly-varying MA($\infty$) arbitrarily close to any valid TPDF.
Using the polar geometry of regular variation, we develop prediction intervals when predicted values are large.


There are many avenues for future work.
A prevalent method for identifying the orders of AR and MA in classical time series is through looking at autocorrelation function (ACF) and partial autocorrelation function (PACF) plots. This necessitates the development of a PACF analogue for our models. Some initial results on partial tail correlation in the finite dimensions have been worked on by Lee \& Cooley in their upcoming paper.

There might also be motivation in extremes to think about a non-causal time series. 
The innovations algorithm looks only backward in time, by definition.
In the simulated time series in our causal setting, we see that extreme behavior is characterized by an asymmetric spike, where large values have subsequent large values, but no large values precede the spike.
Although these transformed linear values have shown good ability to capture dependence summary measures, behavior of the simulated series differs from that of the environmental time series we have so far explored.
A relaxation of causality could likely address this.


\section*{Acknowledgement}
Nehali Mhatre and Daniel Cooley were both partially supported by National Science Foundation award DMS-1811657.

\section*{Supplementary Material}
\label{SM}

Supplementary materials demonstrate that $\mathbb{V}$ is a vector space, that the defined inner product meets the considitions of inner product, and that $T$ is an linear map creatig the isomorphism between $\mathbb{V}$ and $\ell_1$.
The windspeed anomalies data and the R codes for the simulation study results (Section \ref{simulation_study}) and the application to the windspeed anomalies data (Section \ref{sec:app_prediction}) are available at (\href{https://www.stat.colostate.edu/~cooleyd/TransLinTS/}{\textit{https://www.stat.colostate.edu/~cooleyd/TransLinTS/}}).

\bibliographystyle{agsm}
\bibliography{References.bib}

\end{document}